\documentclass[a4paper]{amsart}
\usepackage{amssymb}
\newtheorem{theorem}{Theorem}[section]
\newtheorem{lemma}[theorem]{Lemma}
\newtheorem{corollary}[theorem]{Corollary}
\newtheorem{example}[theorem]{Example}
\newtheorem{definition}[theorem]{Definition}
\title[On the Intuitionistic fuzzy topological spaces]{ on the Intuitionistic fuzzy topological (metric and normed) spaces}
\author[ Reza Saadati ]{Reza Saadati }
\thanks{} 
\address{Department of Mathematics, Azad University, Amol, P.O.Box 678,
Iran}
 \email{rsaadati@eml.cc}

\begin{document}
\maketitle
\begin{abstract}
In this paper, we define precompact set in intuitionistic fuzzy
metric spaces and prove that any subset  of an intuitionistic
fuzzy metric space  is compact if and only if it is precompact and
complete. Also we define topologically complete intuitionistic
fuzzy metrizable spaces and prove that any $G_{\delta }$ set in a
complete intuitionistic fuzzy metric spaces is a topologically
complete intuitionistic fuzzy metrizable space and vice versa.
Finally, we define intuitionistic fuzzy normed spaces and fuzzy
boundedness for linear operators and so we prove that every finite
dimensional intuitionistic fuzzy normed space is complete.
\end{abstract}

\section{Preliminaries}
The theory of fuzzy sets was introduced by L.Zadeh in 1965
\cite{za}. After the pioneering work of Zadeh, there has been a
great effort to obtain fuzzy analogues of classical theories.
Among other fields, a progressive developments is made in the
field of fuzzy topology. The concept of fuzzy topology may have
very important applications in quantum particle physics
particularly in connections with both string and
$\epsilon^{(\infty)}$ theory which were given and studied by
Elnaschie \cite{el1} and \cite{el2}. One of the most important
problems in fuzzy topology is to obtain an appropriate concept of
intuitionistic fuzzy metric space. This problem has been
investigated by J. H. Park \cite{pa}. He has introduced and
studied a notion of intuitionistic fuzzy metric space. We recall
it.

\begin{definition}A binary operation
$*:[0,1]\times[0,1]\longrightarrow[0,1]$ is a continuous t-norm if
it satisfies the following conditions

(1)\, $*$ is associative and commutative,

(2)\, $*$ is continuous,

(3)\, $a*1=a$ for all $a\in[0,1],$

(4)\, $a*b\leq c*d$ whenever $a\leq c$ and $b\leq d,$ for each
$a,b,c,d\in[0,1]$.

\end{definition}

\begin{example}
Two typical examples of continuous t-norm are $a*b=ab$ and
 $a*b=\min(a,b)$.
\end{example}
\begin{definition}
A binary operation $\diamond:[0,1]\times[0,1]\longrightarrow[0,1]$
is a continuous t-conorm if it satisfies the following conditions

(1)\, $\diamond$ is associative and commutative,

(2)\, $\diamond$ is continuous,

(3)\, $a\diamond 0=a$ for all $a\in[0,1],$

(4)\, $a\diamond b\leq c\diamond d$ whenever $a\leq c$ and $b\leq
d,$ for each $a,b,c,d\in[0,1]$.

\end{definition}
\begin{example}
Two typical examples of continuous t-conorm are $a\diamond
b=\min(a+b,1)$ and
 $a\diamond b=\max(a,b)$.
\end{example}
\begin{lemma} If $\ast$ is a continuous t-norm,  $\diamond$ is a continuous
t-conorm and
 $r_{i}\in (0,1), \, 1\leq i\leq 7$, then\\
(i) If $r_{1}>r_{2}$, there are $r_{3},r_{4}\in(0,1)$ such that
$r_{1}\ast r_{3}\geq r_{2}$ and $r_{1}\geq r_{2}\diamond r_{4}$.\\
(ii) If $ r_{5}\in(0,1)$, there are $r_{6},r_{7}\in(0,1)$ such
that $r_{6}\ast r_{6}\geq r_5$ and $r_{5}\geq r_{7}\diamond
r_{7}$.
\end{lemma}

The concept of intuitionistic fuzzy metric space is defined by J.
H. Park \cite{pa}.
\begin{definition}
A 5-tuple $(X,M,N,*,\diamond)$ is called a intuitionistic fuzzy
metric space if $X$ is an arbitrary (non-empty) set, $*$ is a
continuous t-norm, $\diamond$ a continuous t-conorm and $M,N$ are
fuzzy sets on $X^{2}\times (0,\infty)$, satisfying the following
conditions for each $x,y,z\in X$ and $t,s>0$,

(a)\,$M(x,y,t)+N(x,y,t)\leq 1$,

(b)\, $M(x,y,t)>0$,

(c)\, $M(x,y,t)=1$ if and only if $x=y$,

(d)\, $M(x,y,t)=M(y,x,t)$,

(e)\, $M(x,y,t)*M(y,z,s)\leq M(x,z,t+s)$,

(f)\, $M(x,y,.):(0,\infty)\longrightarrow [0,1]$ is continuous.

(g)\, $N(x,y,t)>0$,

(h)\, $N(x,y,t)=0$ if and only if $x=y$,

(i)\, $N(x,y,t)=N(y,x,t)$,

(j)\, $N(x,y,t)\diamond N(y,z,s)\geq N(x,z,t+s)$,

(k)\, $N(x,y,.):(0,\infty)\longrightarrow [0,1]$ is continuous.
\end{definition}
Then $(M,N)$ is called an intuitionistic fuzzy metric on X. The
functions $M(x, y, t)$ and $N(x, y, t)$ denote the degree of
nearness and the degree of non-nearness between x and y with
respect to t, respectively.\\
 Every fuzzy metric space $(X,M,*)$ is an intuitionistic fuzzy
metric space of the form $(X,M,1-M,*,\diamond)$ such that t-norm
$*$ and t-conorm $\diamond $ are associated \cite{lo}, i.e. $x
\diamond y=1-[(1-x)\ast (1-y)]$ for any $x, y\in X$.\\
 In intuitionistic fuzzy metric space $M(x, y, .)$ is
non-decreasing and $N(x, y, .)$is non-increasing for all $x, y\in
X$.\\
Let $(X,d)$ be a metric space. Denote $a*b=ab$ and $a\diamond
b=\min(a+b,1)$ for all $a,b\in[0,1]$ and let $M_d$ and $N_d$ be
fuzzy sets on $X^{2}\times (0,\infty)$ defined as follows:
\begin{eqnarray*}M_{d}(x,y,t)=\frac{ht^{n}}{ht^{n}
+md(x,y)}&,& N_{d}(x,y,t)=\frac{d(x,y)}{kt^{n}
+md(x,y)},\end{eqnarray*} for all $h,k,m,n\in \mathbb{R}^+$.
 Then $(X,M,N,*,\diamond)$ is an intuitionistic fuzzy metric space.

Let $(X,M,N,*,\diamond)$ be a intuitionistic fuzzy metric space.
For $t>0$, the open ball $B(x,r,t)$ with center $x\in X$ and
radius $0<r<1$ is defined by
$$B(x,r,t)=\{y\in X:M(x,y,t)>1-r, N(x,y,t)<r\}.$$

Let $(X,M,N,*,\diamond)$ be a intuitionistic fuzzy metric space.
Let $\tau_{(M,N)}$ be the set of all $A\subset X$ with $x\in A$ if
and only if there exist $t>0$ and $0<r<1$ such that
$B(x,r,t)\subset A$. Then $\tau_{(M,N)}$ is a topology on $X$
(induced by the intuitionistic fuzzy metric $(M,N))$. This
topology is Hausdorff and first countable. A sequence $\{x_{n}\}$
in $X$ converges to $x$ if and only if $M(x_{n},x,t)\to 1$ and
$N(x_{n},x,t)\to 0$ as $n\to\infty$, for each $t>0$. It is called
a Cauchy sequence if for each $0<\varepsilon<1$ and $t>0$, there
exits $n_{0}\in \mathbb{N}$ such that
$M(x_{n},x_{m},t)>1-\varepsilon$ and
$N(x_{n},x_{m},t)<\varepsilon$ for each $n,m\geq n_{0}$. The
intuitionistic fuzzy metric space $(X,M,N,*,\diamond)$ is said to
be complete if every Cauchy sequence is convergent. A subset $A$
of $X$ is said to be IF-bounded if there exists $t>0$ and $0<r<1$
such that $M(x,y,t)>1-r$  and $N(x,y,t)<r$ for all $x,y \in A$.

A collection $\mathcal{U}$ of open sets is called an open cover of
$A$ if $A\subseteq \bigcup_{U\in \mathcal{U}} U$. A subspace $A$
of an intuitionistic fuzzy metric space $(X,M,N,*,\diamond)$ is
compact if every open cover of $A$ has a finite subcover. If every
sequence in $A$ has a convergent subsequence to a point in $A$
then it is called \textit{sequential compact}.

\begin{theorem}
\cite{pa}. In a intuitionistic fuzzy metric space every compact
set is closed and IF-bounded.
\end{theorem}
\begin{corollary}\cite{pa}.
Every closed subset of a complete intuitionistic fuzzy metric
space is complete.
\end{corollary}

\section{Precompact Intuitionistic Fuzzy Metric Spaces}
\begin{definition}
Let $(X,M,N,*,\diamond)$ be a intuitionistic fuzzy metric space
and $A\subset X$. We say $A$ is precompact if for each $0<r<1$ and
$t>0$ there exists a finite subset $S$ of $A$ such that $$
A\subseteq \bigcup_{x\in S} B(x,r,t).$$
\end{definition}
\begin{lemma} Let $(X,M,N,*,\diamond)$ be a intuitionistic fuzzy metric space
and $A\subset X$. $A$ is a precompact set if and only if for every
$0<r<1$ and $t>0$, there exists a finite subset $S$ of $X$ such
that
\begin{eqnarray} A\subseteq \bigcup_{x\in
S}B(x,r,t).\end{eqnarray}\end{lemma}

\textbf{Proof.} Let $0<r<1$ and $t>0$ and condition (1) holds. By
continuity of $\ast,\diamond$, there exists $s\in(0,1)$ such that
$(1-s)*(1-s)> 1-r$ and $s\diamond s< r$ . Now we applying
condition (1) for $s$ and $\frac{t}{2}$, there exists a subset
$S'=\{x_{1},...,x_{n}\}$ of $X$ such that $A\subseteq
\bigcup_{x_{i}\in S'}B(x_i,s,\frac{t}{2})$. We assume that
$B(x_j,s,\frac{t}{2})\cap A\neq \phi$, otherwise we omit $x_{j}$
from $S'$ and so we have $A\subseteq \bigcup_{x_{i}\in
S'-\{x_{j}\}}B(x_i,s,\frac{t}{2})$. For every $i=1,...,n$ we
select $y_{i}$ in $B(x_i,s,\frac{t}{2})\cap A$, and we put
$S=\{y_{1},...,y_{n}\}$. Now for every $y$ in $A$, there exists
$i\in \{1,...,n\}$ such that $M(y,x_i,\frac{t}{2})>1-s$ and
$N(y,x_i,\frac{t}{2})<s$ . Therefore we have
\begin{eqnarray*}M(y,y_i,t)&>&
M(y,x_i,\frac{t}{2})*M(x_i,y_i,\frac{t}{2})\\
&>& (1-s)\ast (1-s)\\
&>&1-r,\end{eqnarray*} and
\begin{eqnarray*}N(y,y_i,t)&<&
N(y,x_i,\frac{t}{2})\diamond N(x_i,y_i,\frac{t}{2})\\
&<& s\diamond s\\
&<&r.\end{eqnarray*}
 Which implies that $A\subseteq
\bigcup_{p_{i}\in S}B(x,r,t)$. The converse is trivial.\qed

\begin{lemma} Let $(X,M,N,*,\diamond)$ be a intuitionistic fuzzy metric space
and $A\subset X$. If $A$ is a precompact set then so is its
closure \={A}.\end{lemma}

\textbf{Proof.} Let $r\in(0,1)$ and $t>0$, then by continuity of
$\ast,\diamond$ there exists $s\in(0,1)$ such that $(1-s)*(1-s)>
1-r$ and $s\diamond s< r$, also there exists a finite subset
$S'=\{x_{1},...,x_{n}\}$ of $X$  such that $A\subseteq
\bigcup_{x_{i}\in S'}B(x_i,s,t/2).$ But for every $y$ in \={A}
there exists $x\in A$ such that $M(x,y,t/2)>1-s$ and
$N(x,y,t/2)<s$ and there exists $1\leq i\leq n$ such that
$M(x,x_i,t/2)>1-s$ and $N(x,x_i,t/2)<s$, therefore
\begin{eqnarray*} M(y,x_i,t)&>&
M(y,x,t/2)*M(x,x_i,t/2)\\
&>&(1-s)*(1-s)\\
&>& 1-r,\end{eqnarray*} and
\begin{eqnarray*} N(y,x_i,t)&<&
N(y,x,t/2)\diamond N(x,x_i,t/2)\\
&<&s\diamond s\\
&<& r.\end{eqnarray*}
 Hence \={A}$\subseteq \bigcup_{x_{i}\in
S}B(x_i,r,t)$, i.e. \={A} is precompact set.\qed
\begin{theorem} Let $(X,M,N,*,\diamond)$ be a intuitionistic fuzzy metric space
and $A\subset X$. $A$ is a precompact set if and only if every
sequence has a  Cauchy subsequence.\end{theorem}
 \textbf{Proof.} Let $A$ be a precompact set. Let $\{p_{n}\}$ be a sequence in $A$.
 For every $k\in \mathbb{N}$, there
exists a finite subset $S_{k}$ of $X$ such that $ A\subseteq
\bigcup_{x\in S_{k}}B(x,\frac{1}{k},\frac{1}{k})$. Hence, for
$k=1$, there exists $x_{1}\in S_{1}$ and a subsequence
$\{p_{1,n}\}$ of $\{p_{n}\}$ such that $p_{1,n}\in B(x_1,1,1)$,
for every $n\in \mathbb{N}$. Similarly, there exists $x_{2}\in
S_{2}$ and a subsequence $\{p_{2,n}\}$ of $\{p_{1,n}\}$ such that
$p_{2,n}\in B(x_2,\frac{1}{2},\frac{1}{2})$, for every $n\in
\mathbb{N}$. Continuing this process, we get $x_{k}\in S_{k}$ and
subsequences $\{p_{k,n}\}$ of $\{p_{k-1,n}\}$ such that
$p_{k,n}\in B(x,\frac{1}{k},\frac{1}{k})$, for every $n\in
\mathbb{N}$. Now we consider the subsequence $\{p_{n,n}\}$ of
$\{p_{n}\}$. For every $r\in(0,1)$ and $t>0$, by continuity of
$\ast,\diamond$, there exists an $n_{0}\in\mathbb{N}$ such that
$(1-\frac{1}{n_{0}})*(1-\frac{1}{n_{0}})>1-r$,
$\frac{1}{n_{0}}\diamond \frac{1}{n_{0}}< r$ and
$\frac{2}{n_{0}}<t$. Therefore for every $l,m\geq n_{0}$, we have
\begin{eqnarray*} M(p_{l,l},p_{m,m},t)&\geq &
M(p_{l,l},p_{m,m},\frac{2}{n_{0}})\\
&\geq & M(p_{l,l},x_{n_{0}},\frac{1}{n_{0}})\ast M(x_{n_{0}},p_{m,m},\frac{1}{n_{0}})\\
&>& (1-\frac{1}{n_{0}})*(1-\frac{1}{n_{0}})\\
&> & 1-r,\end{eqnarray*}and
\begin{eqnarray*} N(p_{l,l},p_{m,m},t)&\leq &
N(p_{l,l},p_{m,m},\frac{2}{n_{0}})\\
&\leq & N(p_{l,l},x_{n_{0}},\frac{1}{n_{0}})\ast N(x_{n_{0}},p_{m,m},\frac{1}{n_{0}})\\
&<& \frac{1}{n_{0}}\diamond\frac{1}{n_{0}}\\
&< & r.\end{eqnarray*}
 Hence $\{p_{n,n}\}$ is a Cauchy sequence in $(X,M,N,*,\diamond)$.

Conversely, suppose that $A$ is not a precompact set. Then there
exists $r\in(0,1)$ and $t>0$ such that for every finite subset $S$
of $X$, $A$ is not a subset of $\bigcup_{x\in S}B(x,r,t)$. Fix
$p_{1}\in A$. Since $A$ is not a subset of $ \bigcup_{x\in
\{p_{1}\}}B(x,r,t)$, there exists $p_{2}\in A$ such that
$M(p_{1},p_{2},t)\leq 1-r$ and $N(p_{1},p_{2},t)\geq r$. Since $A$
is not a subset of $ \bigcup_{x\in\{p_{1},p_{2}\}}B(x,r,t)$, there
exists a $p_{3}\in A$ such that $M(p_{1},p_{3},t)\leq 1-r$,
$N(p_{1},p_{3},t)\geq r$ and $M(p_{2},p_{3},t)\leq 1-r$,
$N(p_{2},p_{3},t)\geq r$. Continuing this process, we construct a
sequence $\{p_{n}\}$ of distinct points in $A$ such that
$M(p_{i},p_{j},t)\leq 1-r$ and $N(p_{i},p_{j},t)\geq r$ for every
$i\neq j$. Therefore $\{p_{n}\}$ has
not Cauchy subsequence.\qed\\
\begin{lemma}Let $(X,M,N,*,\diamond)$ be a intuitionistic fuzzy metric space.
If a Cauchy sequence clusters to a point $x\in X$, then the
sequence converges to $x$.\end{lemma} \textbf{Proof.} Let
$\{x_n\}$ be a Cauchy sequence in $(X,M,N,*,\diamond)$ having a
cluster point $x\in X$. Then, there is a subsequence
$\{x_{n_{k}}\}$ of $\{x_n\}$ that converges to $x$ with respect to
$\tau_{(M,N)}$. Thus, given $r$, with $0<r<1$ and $t>0$ there is
an $N\in \mathbb{N}$ such that for each $k\geq N$,
$M(x,x_{n_{k}},t/2)>1-s$ and $N(x,x_{n_{k}},t/2)<s$ where
$s\in(0,1)$ and satisfies $(1-s)*(1-s)> 1-r$ and $s\diamond s< r$.
 On the other hand, there is $n_{1}\geq n_{N}$ such that for
 each $n,m\geq n_{1}$, we have $M(x_m,x_n,t/2)>1-s$ and $N(x_m,x_n,t/2)<s$. Therefore,
 for each $n\geq n_1$, we have
 \begin{eqnarray*} M(x,x_n,t)& \geq
 & M(x,x_{n_{k}},t)* M(x_{n_{k}},x_n,t)\\
 &>&(1-s)*(1-s)\\
&>& 1-r,\end{eqnarray*} and
\begin{eqnarray*} N(x,x_n,t)& \leq
 & N(x,x_{n_{k}},t)\diamond N(x_{n_{k}},x_n,t)\\
 &<&s\diamond s\\
&<& r.\end{eqnarray*} We conclude that the Cauchy sequence
$\{x_n\}$ converges to $x$.\qed
\begin{lemma}Let $(X,M,N,*,\diamond)$ be a intuitionistic fuzzy metric space.
Then $(X,\tau_{(M,N)})$ is a metrizable topological
space.\end{lemma} \textbf{Proof.} For each $n\in \mathbb{N}$
define \begin{eqnarray*} U_n= \{(x,y)\in X\times X:
&M(x,y,\frac{1}{n})>1-\frac{1}{n}
,&N(x,y,\frac{1}{n})<\frac{1}{n}\}.\end{eqnarray*} We sall prove
that $\{U_n:n\in \mathbb{N}\}$ is a base for a uniformity
$\mathcal{U}$ on $X$ whose induced topology coincides with
$\tau_{(M,N)}$. We first note that for each $n\in \mathbb{N}$,
$\{(x,x):x\in X\}\subseteq U_n$, $U_{n+1}\subseteq U_n$ and
$U_n=U_{n}^{-1}$.

On the other hand, for each $n\in \mathbb{N}$, there is, by the
continuity of $\ast ,\diamond$, an $m\in \mathbb{N}$ such that
$m>2n$, $(1-\frac{1}{m})*(1-\frac{1}{m})>1-\frac{1}{n}$ and
$\frac{1}{m}\diamond\frac{1}{m}<\frac{1}{n}$. Then, $U_m \circ
U_m\subseteq U_n$. Indeed, let $(x,y)\in U_m$ and $(y,z)\in U_m$.
Since $M(x,y,.)$ and $N(x,y,.)$ are nondecreasing and
nonincreasing, respectively, $M(x,y,\frac{1}{n})\geq
M(x,z,\frac{2}{m})$ and $N(x,y,\frac{1}{n})\leq
N(x,z,\frac{2}{m})$. So \begin{eqnarray*} M(x,z,\frac{1}{n})&\geq
&M(x,y,\frac{1}{m})*M(y,z,\frac{1}{m})\\
&>& (1-\frac{1}{m})*(1-\frac{1}{m})>1-\frac{1}{n},\end{eqnarray*}
and \begin{eqnarray*} N(x,z,\frac{1}{n})&\leq
&N(x,y,\frac{1}{m})\diamond N(y,z,\frac{1}{m})\\
&<& \frac{1}{m}\diamond\frac{1}{m}<\frac{1}{n}.\end{eqnarray*}
Therefore $(x,z) \in U_n$. Thus $\{U_n:n\in \mathbb{N}\}$ is a
base for a uniformity $\mathcal{U}$ on $X$. Since for each $x\in
X$ and each $n\in \mathbb{N}$, \begin{eqnarray*}U_n(x)= \{y\in X:
&M(x,y,\frac{1}{n})>1-\frac{1}{n}
,&N(x,y,\frac{1}{n})<\frac{1}{n}\}=B(x,\frac{1}{n},\frac{1}{n}),\end{eqnarray*}
we deduce that the topology induced by $\mathcal{U}$ coincides
with $\tau_{(M,N)}$. Then $(X,\tau_{(M,N)})$ is a metrizable
topological space.\qed

Note that, in every metrizable space every sequentially copmact
set is compact.
\begin{corollary} A subset $A$ of  intuitionistic fuzzy metric space $(X,M,N,*,\diamond)$ is compact
 if and only if it is precompact and complete.\end{corollary}
\section{Complete Intuitionistic Fuzzy Metric Spaces}
\begin{lemma}
Let $(X,M,N,*,\diamond)$ be a intuitionistic fuzzy metric space
and let $\lambda,\eta\in(0,1)$ such that $\lambda+\eta\leq 1$ then
there exists a intuitionistic fuzzy metric $(m,n)$ on $X$ such
that $m(x,y,t)\geq\lambda$ and $n(x,y,t)\leq \eta$ for each
$x,y\in X$ and $t>0$ and $(m,n)$ and $(M,N)$ induce the same
topology on $X$.
\end{lemma}
\textbf{Proof.} We define $m(x,y,t)=\max\{\lambda,M(x,y,t)\}$ and
$n(x,y,t)=\min\{\eta,N(x,y,t)\}$. We claim that $(m,n)$ is
intuitionistic fuzzy metric on $X$. The properties of (a),(b),(c),
(d),(f),(g),(h),(i) and (k) are immediate from the definition. For
triangle inequalities, suppose that $x,y,z\in X $ and $t,s>0$.
Then $ m(x,z,t+s)\geq \lambda$ and so $m(x,z,t+s)\geq
m(x,y,t)*m(y,z,s)$ when either $m(x,y,t)=\lambda$ or
$m(y,z,s)=\lambda$. The only remaining case is when
$m(x,y,t)=M(x,y,t)>\lambda$ and $m(y,z,s)=M(y,z,s)>\lambda$. But
$M(x,z,t+s)\geq M(x,y,t)*M(y,z,s)$ and $m(x,z,t+s)\geq M(x,z,t+s)$
and so $m(x,z,t+s)\geq m(x,y,t)*m(y,z,s)$. Also, then $
n(x,z,t+s)\leq \eta$ and so $n(x,z,t+s)\leq n(x,y,t)\diamond
n(y,z,s)$ when either $n(x,y,t)=\eta$ or $n(y,z,s)=\eta$. The only
remaining case is when $n(x,y,t)=N(x,y,t)<\eta$ and
$n(y,z,s)=N(y,z,s)<\eta$. But $N(x,z,t+s)\leq N(x,y,t)\diamond
N(y,z,s)$ and $n(x,z,t+s)\leq N(x,z,t+s)$ and so $n(x,z,t+s)\leq
n(x,y,t)\diamond n(y,z,s)$. Thus $(m,n)$ is a intuitionistic fuzzy
metric on $X$. It only remains to show that the topology induced
by $(m,n)$ is the same as that induced by $(M,N)$. But we have $
m(x_{n},x,t)\longrightarrow 1$ and $n(x_{n},x,t)\longrightarrow 0$
if and only if $\{\lambda,M(x_{n},x,t)\}\longrightarrow 1$ and
$\{\eta,N(x_{n},x,t)\}\longrightarrow 0$ if and only if
$M(x_{n},x,t)\longrightarrow 1$ and $N(x_{n},x,t)\longrightarrow
0$, for each $t>0$, and we are done.\qed

The intuitionistic fuzzy metric $(m,n)$ in above lemma is said to
be bounded by $(\lambda,\eta)$.

\begin{definition}
Let $(X,M,N,*,\diamond)$ be a intuitionistic fuzzy metric space,
$x\in X$ and $\phi\neq A\subseteq X$. We define

$$D(x,A,t)=\sup\{M(x,y,t):y\in A\}\quad(t>0),$$ and
$$C(x,A,t)=\inf\{N(x,y,t):y\in A\}\quad(t>0).$$
Note that $D(x,A,t)$  and $C(x,A,t)$ are a degree of closeness and
non closeness of $x$ to $A$ at $t$, respectively.
\end{definition}

\begin{definition}
A topological space is called a topologically complete
intuitionistic fuzzy metrizable space if there exists a
 complete intuitionistic fuzzy metric inducing the
given topology on it.
\end{definition}
\begin{example} Let $X=(0,1]$. The intuitionistic fuzzy metric space $(X,M,N,\min,\max)$
where $M(x,y,t)=\frac{t}{t+|x-y|}$ and
$N(x,y,t)=\frac{|x-y|}{t+|x-y|}$ (standard intuitionistic fuzzy
metric, see \cite{pa}) is not complete, because the Cauchy
sequence $\{1/n\}$ in this space is not convergent. Now consider
the 5-tuple $(X,m,n,\min,\max)$, where
$m(x,y,t)=\frac{t}{t+|x-y|+|\frac{1}{x}-\frac{1}{y}|}$ and
$n(x,y,t)=\frac{|x-y|+|\frac{1}{x}-\frac{1}{y}|}{t+|x-y|+|\frac{1}{x}-\frac{1}{y}|}$.
It is straightforward to show that $(X,m,n,\min,\max)$ is a
intuitionistic fuzzy metric space which is complete. Since,
$x_{n}$ tend to $x$ with respect to intuitionistic fuzzy metric
$(M,N)$, if and only if $|x_{n}-x|\longrightarrow 0$, if and only
if $x_{n}$ tends to $x$ with respect to intuitionistic fuzzy
metric $(m,n)$, hence $(M,N)$ and $(m,n)$ are equivalent
intuitionistic fuzzy metrics. Therefore the intuitionistic fuzzy
metric space $(X,M,N,\min,\max)$ is topologically complete
intuitionistic fuzzy metrizable.
\end{example}

\begin{lemma}Intuitionistic fuzzy  metrizability is preserved under countable Cartesian product.
\end{lemma}
\textbf{Proof.} Without loss of generality we may assume that the
index set is $\mathbb N$. Let $
\{(X_{n},m_{n},n_{n},*,\diamond):n\in\mathbb N \}$ be a collection
of intuitionistic fuzzy metrizable spaces. Let $\tau_{n}$ be the
topology induced by $(m_{n},n_{n})$ on $X_{n}$ for $n\in \mathbb
N$ and let $(X,\tau)$ be the Cartesian product of
$\{(X_{n},\tau_{n}):n\in\mathbb N \}$ with product topology. We
have to prove that there is a intuitionistic fuzzy metric $(m,n)$
on $X$ which induces the topology $\tau $. By the above lemma, we
may suppose that $(m_{n},n_{n})$ is bounded by
$(1-\varepsilon^{(n)},\varepsilon^{[n]})$ where
$\varepsilon^{(n)}=\overbrace{\varepsilon*
\varepsilon*\cdots*\varepsilon}^{n}$,
$\varepsilon^{[n]}=\overbrace{\varepsilon\diamond
\varepsilon\diamond\cdots\diamond\varepsilon}^{n}$ and
$\varepsilon\in(0,1)$ (see, \cite{as2}), i.e.
$m_{n}(x_{n},y_{n},t)=max\{1-\varepsilon^{(n)},M(x_{n},y_{n},t)\}$
and
$n_{n}(x_{n},y_{n},t)=\min\{\varepsilon^{[n]},N(x_{n},y_{n},t)\}$.
Points of $X=\prod_{n\in N}X_{n}$ are denoted as sequences
$x=\{x_{n} \}$ with $x_{n}\in X_{n}$ for $n\in\mathbb N$. Define
$m(x,y,t)=\prod_{n=1}^{\infty}m_{n}(x_{n},y_{n},t)$ and
$n(x,y,t)=\coprod_{n=1}^{\infty}n_{n}(x_{n},y_{n},t)$ , for each
$x,y\in X$ and $t>0$ where $\prod_{n=1}^{m}
a_{n}=a_{1}*a_{2}*\cdots*a_{m}$ and $\coprod_{n=1}^{m}
a_{n}=a_{1}\diamond a_{2}\diamond\cdots\diamond a_{m}$. First note
that $(m,n)$ is well defined since
$a_{i}=\prod_{n=1}^{i}(1-\varepsilon^{(n)})$ is decreasing and
bounded then converges to $\alpha\in(0,1)$ also
$b_{i}=\coprod_{n=1}^{i}\varepsilon^{[n]}$ is increasing and
bounded then converges to $\beta\in(0,1)$. Also $(m,n)$ is a
intuitionistic fuzzy metric on $X$ because each $(m_{n},n_{n})$ is
a intuitionistic fuzzy metric. Let $\mathcal{U}$ be the topology
induced by intuitionistic fuzzy metric $(m,n)$. We claim that
$\mathcal{U}$ coincides with $\tau$. If $G\in \mathcal{U}$ and
$x=\{x_{n}\}\in G$, then there exists $0<r<1$ and $t>0$ such that
$B(x,r,t)\subset G$. For each $0<r<1$, we can find a sequence
$\{\delta_{n}\}$ in $(0,1)$ and a positive integer $N_{0}$ such
that

$$\prod_{n=1}^{N_{0}} (1-\delta_{n})*\prod_{n=N_{0}+1}^{\infty}
(1-\varepsilon^{(n)})>1-r,$$ and
$$\coprod_{n=1}^{N_{0}} \delta_{n}\diamond\coprod_{n=N_{0}+1}^{\infty}
\varepsilon^{[n]}<r.$$

For each $n=1,2,\cdots,N_{0}$, let $V_{n}=B(x_{n},\delta_{n},t)$,
 where the ball is with respect to intuitionistic fuzzy metric $(m_{n},n_{n})$. Let
$V_{n}=X_{n}$ for $n>N_{0}$. Put $V=\prod_{n\in N}V_{n}$, then
$x\in V$ and $V$ is an open set in the product topology $\tau$ on
$X$. Furthermore $V\subset B(x,r,t)$, since for each $y\in V$
\begin{eqnarray*}
m(x,y,t)&=&\prod_{n=1}^{\infty}
m_{n}(x_{n},y_{n},t)\\
&=&\prod_{n=1}^{N_{0}}m_{n}(x_{n},y_{n},t)*\prod_{n=N_{0}+1}^{\infty} m_{n}(x_{n},y_{n},t)\\
& \geq & \prod_{n=1}^{N_{0}}
(1-\delta_{n})*\prod_{n=N_{0}+1}^{\infty}
(1-\varepsilon^{(n)})\\
&>& 1-r,
\end{eqnarray*}
and
\begin{eqnarray*}
n(x,y,t)&=&\coprod_{n=1}^{\infty}
n_{n}(x_{n},y_{n},t)\\
&=&\coprod_{n=1}^{N_{0}}n_{n}(x_{n},y_{n},t)\diamond\coprod_{n=N_{0}+1}^{\infty} n_{n}(x_{n},y_{n},t)\\
& \leq & \coprod_{n=1}^{N_{0}}
\delta_{n}\diamond\coprod_{n=N_{0}+1}^{\infty}
\varepsilon^{[n]}\\
&<& r.
\end{eqnarray*}
 Hence $V\subset B(x,r,t)\subset G$. Therefore $G$ is open in
the product topology.

Conversely, suppose $G$ is open in the product topology and let
$x=\{x_{n}\}\in G$. Choose a standard basic open set $V$ such that
$x\in V$ and $V\subset G$. Let $V=\prod_{n\in N}V_{n}$, where each
$V_{n}$ is open in $X_{n}$ and $V_{n}=X_{n}$ for all $n>N_{0}$.
For $n=1,2,\cdots,N_{0}$, let $1-r_{n}=D_{n}(x_{n},X_{n}-V_{n},t)$
and $q_{n}=C_{n}(x_{n},X_{n}-V_{n},t)$, if $X_{n}\neq V_{n}$, and
$r_{n}=\varepsilon^{(n)}$ and $q_{n}=\varepsilon^{[n]}$,
otherwise. Let $r=\min\{r_{1},r_{2},\cdots,r_{N_{0}}\}$ ,
$q=\min\{q_{1},q_{2},\cdots,q_{N_{0}}\}$ and $p=\min\{r,q\}$. We
claim that $B(x,p,t)\subset V$. If $y=\{y_{n}\}\in B(x,p,t)$, then
$m(x,y,t)=\prod_{n=1}^{\infty} m_{n}(x_{n},y_{n},t)>1-p$   and so
$ m_{n}(x_{n},y_{n},t)>1-p\geq 1-r\geq 1-r_{n}$ and
$n(x,y,t)=\coprod_{n=1}^{\infty} n_{n}(x_{n},y_{n},t)<p$   and so
$ n_{n}(x_{n},y_{n},t)<p\leq q\leq q_{n}$ for each
$n=1,2,\cdots,N_{0}$. Then $y_{n}\in V_{n}$, for
$n=1,2,\cdots,N_{0}$. Also for $n>N_{0}, y_{n}\in V_n=X_{n}$.
Hence $y\in V$ and so $B(x,p,t)\subset V \subset G$. Therefore $G$
is open with respect to the intuitionistic fuzzy metric topology
and $\tau
 \subset \mathcal{U}$. Hence $\tau$ and $\mathcal{U}$
 coincide.\qed

\begin{theorem}
An open subspace of a complete intuitionistic fuzzy metrizable
space is a topologically complete intuitionistic fuzzy metrizable
space.
\end{theorem}
\textbf{Proof.} Let $(X,M,N,*,\diamond)$ be a complete
intuitionistic fuzzy metric space and $G$ an open subspace of $X$.
If the restriction of $(M,N)$ to $G$ is not complete we can
replace $(M,N)$ on $G$ by other intuitionistic fuzzy metric as
follows. Define $f:G\times (0,\infty)\longrightarrow R^{+}$ by
$f(x,t)=\frac{1}{1-D(x,X-G,t)}$  ($f$ is undefined if $X-G$ is
empty, but then there is nothing to prove.) Fix an arbitrary $s>0$
 and for $x,y\in G$ define

$$m(x,y,t)=M(x,y,t)*M(f(x,s),f(y,s),t),$$ and
$$n(x,y,t)=N(x,y,t)$$
for each $t>0$. We claim that $(m,n)$ is intuitionistic fuzzy
metric on $G$. The properties (a),(b),(c),(d),(f),(g),(h),(i),(j)
and (k) are immediate from the definition.
 For triangle inequality (e), suppose that $x,y,z \in G$ and
 $t,s,u>0$, then
\begin{align*}
m(x,y,t)&*m(y,z,u)=\\
&=(M(x,y,t)*M(f(x,s),f(y,s),t))*(M(y,z,u)*M(f(y,s),f(z,s),u))\\
&=(M(x,y,t)*M(y,z,u))*(M(f(x,s),f(y,s),t)*M(f(y,s),f(z,s),u))\\
&\leq M(x,z,t+u)*M(f(x,s),f(z,s),t+u) =m(x,z,t+u).
\end{align*}

We show that $(m,n)$ and $(M,N)$ are equivalent intuitionistic
fuzzy metrics on $G$. We do this by showing that
$m(x_{n},x,t)\longrightarrow 1$ if and only if
$M(x_{n},x,t)\longrightarrow 1$ and $n(x_{n},x,t)\longrightarrow
0$ if and only if $N(x_{n},x,t)\longrightarrow 1$ of course the
second part is trivial. Since $m(x,y,t)\leq M(x,y,t)$ for all
$x,y\in G$ and $t>0$, $M(x_{n},x,t)\longrightarrow 1$ whenever
$m(x_{n},x,t)\longrightarrow 1$. To prove the converse, let
$M(x_{n},x,t)\longrightarrow 1$, we know from \cite{lr}
Proposition 1,  $M$ is continuous function on $X\times X\times
(0,\infty)$, then since
\begin{eqnarray*}
\lim_{n} D(x_{n},X-G,s)&=&\lim_{n}(\sup\{M(x_{n},y,s):y\in G\})\\
&\geq &\lim_{n}M(x_{n},y,s)\\
&=&M(x,y,s). \end{eqnarray*} Therefore $\lim_{n}D(x_{n},X-G,s)\geq
D(x,X-G,s)$. On the other hand, there exists a $y_{0}\in X-G$ and
$n_{0}\in \mathbb{N}$ such that for every $n\geq n_{0}$ we have
                   $$  D(x_{n},X-G,s)\ast (1-\frac{1}{n})\leq M(x_{n},y_{0},s).$$
Then $\lim_{n}D(x_{n},X-G,s)\leq M(x,y_{0},s)\leq
\sup\{M(x,y,s):y\in X-G \}=D(x,X-G,s)$. Therefore
$\lim_{n}D(x_{n},X-G,s)=D(x,X-G,s)$. This implies
$M(f(x_{n},s),f(x,s),t)\longrightarrow 1$. Hence
$m(x_{n},x,t)\longrightarrow 1$. Therefore $(m,n)$ and $(M,N)$ are
equivalent. Next we show that $(m,n)$ is a complete intuitionistic
fuzzy metric. Suppose that $\{x_{n}\}$ is a Cauchy sequence in $G$
with respect to $(m,n)$. Since for each $m,n\in\mathbb N$, and
$t>0$ $m(x_{m},x_{n},t)\leq M(x_{m},x_{n},t)$ and
$n(x,y,t)=N(x,y,t)$, therefore $\{x_{n}\}$ is also a Cauchy
sequence with respect to $(M,N)$. By completeness of
$(X,M,N,*,\diamond)$, $\{x_{n}\}$ converges to point  $p$ in $X$.
We claim that $p\in G$. Assume otherwise, then for each
$n\in\mathbb N$, if $p\in X-G$ and $M(x_{n},p,t)\leq
D(x_{n},X-G,t)$, then
$$1-M(x_{n},p,t)\geq 1-D(x_{n},X-G,t)>0,$$
Therefore $$\frac{1}{1-D(x_{n},X-G,t)}\geq
\frac{1}{1-M(x_{n},p,t)},$$  That is
$$f(x_{n},t)\geq\frac{1}{1-M(x_{n},p,t)},$$
for each $t>0$. Therefore as $n\longrightarrow \infty$, for every
$t>0$ we get
 $f(x_{n},t)\longrightarrow \infty$. In particular,
 $f(x_{n},s)\longrightarrow\infty$. On the other hand,
 $M(f(x_{n},s),f(x_{m},s),t)\geq m(x_{m},x_{n},t)$, for every $m,n\in\mathbb
 N$, that is
 $\{f(x_{n},s)\}$ is an F-bounded sequence (see, \cite{gv1}). This contradiction shows that
 $p\in G$. Hence $\{x_{n}\}$ converges to $p$ with respect to $(m,n)$ and $(G,m,n,*,\diamond)$ is a complete
intuitionistic fuzzy metrizable space.\qed

 \begin{corollary}
 A $G_{\delta}$ set in a complete intuitionistic fuzzy metric space is a topologically complete intuitionistic fuzzy
 metrizable space.
 \end{corollary}
 \begin{theorem} Let $(Y,M,N,*,\diamond)$ be a intuitionistic fuzzy metric space and $X$ be a topologically
complete intuitionistic fuzzy metrizable subspace of $Y$. Then $X$
is a $G_{\delta}$ subset of $Y$.\end{theorem} \textbf{Proof.} Let
$(X,M',N',*,\diamond)$ be intuitionistic fuzzy metric space that
induces the same topology for $X$ as does $(M,N)$. For each $x\in
X$ and each $n\in \mathbb{N}$, let $r_{n}(x)$ be a positive real
number such that, $r_{n}(x)<\frac{1}{n}$ and
$M'(w,x,t)>1-\frac{1}{n}$ and $N'(w,x,t)<\frac{1}{n}$, whenever
$w\in X$ and $M(w,x,t)>1-r_{n}(x)$ and $N(w,x,t)<r_{n}(x)$ for
each $t>0$. Suppose that
$\mathcal{G}_{n}=\cup_{n}\{B_{(M,N)}(x,r_{n}(x),t):x\in X,t>0\}$
for each $n\in \mathbb{N}$, and $\Gamma
=\cap_{n}\{\mathcal{G}_{n}:n\in \mathbb{N}\}$. Then $\Gamma$ is a
$G_{\delta}$ subset of $Y$ which clearly contains $X$. It is
enough to shown that $\Gamma \subseteq X$.

Let $x_{0}\in\Gamma$, then $x_{0}\in \mathcal{G}_{n}$ for each
$n\in \mathbb{N}$. Hence for each $n\in \mathbb{N}$, there is
$x_{n}\in X$ such that $x_{0}\in B_{(M,N)}(x_{n},r_{n}(x),t)$.
Therefore $M(x_0,x_n,t)>1-r_{n}(x)>1-\frac{1}{n}$ and
$N(x_0,x_n,t)<r_{n}(x)<\frac{1}{n}$ for each $n\in \mathbb{N}$ and
$t>0$. This means that $x_{n}\longrightarrow x_{0}$ in $Y$.

Now, let $0<\varepsilon<1$ and $N\in \mathbb{N}$ such that
$(1-\frac{1}{N})*(1-\frac{1}{N})>1-\varepsilon$ and
$\frac{1}{N}\diamond\frac{1}{N}<\varepsilon$. Let $m\in
\mathbb{N}$ be such that
$$(1-\frac{1}{m})*M(x_0,x_N,t)>1-r_{N}(x_N),$$ and $$\frac{1}{m}\diamond N(x_0,x_N,t)<r_{N}(x_N),$$
 Now for every $k\in \mathbb{N}$ and $k>m$ we have
\begin{eqnarray*}
M(x_k,x_N,2t)&\geq & M(x_k,x_0,t)*M(x_0,x_N,t)\\
& >& (1-\frac{1}{k})*M(x_0,x_N,t)\\
& \geq & (1-\frac{1}{m})*M(x_0,x_N,t)\\
& > & 1-r_N(x_N),\end{eqnarray*} and
\begin{eqnarray*}
N(x_k,x_N,2t)&\leq & N(x_k,x_0,t)\diamond N(x_0,x_N,t)\\
& <& (\frac{1}{k})\diamond N(x_0,x_N,t)\\
& \leq & (\frac{1}{m})\diamond N(x_0,x_N,t)\\
& < & r_N(x_N),\end{eqnarray*} Therefore
$M'(x_k,x_N,2t)>1-\frac{1}{N}$ and $N'(x_k,x_N,2t)<\frac{1}{N}$ .
If $k,l>m$, then
\begin{eqnarray*}
M'(x_k,x_l,4t)&\geq & M'(x_k,x_N,2t)*M'(x_N,x_l,2t)\\
& >& (1-\frac{1}{N})*(1-\frac{1}{N})>1-\varepsilon,\end{eqnarray*}
and
\begin{eqnarray*}
N'(x_k,x_l,4t)&\leq & N'(x_k,x_N,2t)\diamond N'(x_N,x_l,2t)\\
& <& \frac{1}{N}\diamond\frac{1}{N}< \varepsilon.\end{eqnarray*}
  Hence the sequence $\{x_{n}\}$
is Cauchy in the complete intuitionistic fuzzy metric space
$(X,M',N',*,\diamond)$ and so  convergent to some member of $X$.
Since $x_{n}\longrightarrow x_{0}$ in $Y$, it follows that
$x_{0}\in X$, so $\Gamma \subseteq X$.\qed
\section{INTUITIONISTIC FUZZY NORMED SPACES}
In this section, using the idea of intuitionistic fuzzy metric
space, we define the notion of intuitionistic fuzzy normed spaces
with the help of continuous t-norms and continuous t-conorms as a
generalization of fuzzy normed space due to Saadati and Vaezpour
\cite{sa}.
\begin{definition}The 5-tuple $(V,\mu,\nu,\ast,\diamond)$ is said to be a intuitionistic fuzzy normed space if $V$ is
a vector space, $\ast$ is a continuous t-norm, $\diamond$ is a
continuous t-conorm,  and $\mu,\nu$ are fuzzy sets on
$V\times(0,\infty)$ satisfying the following conditions for every
$x,y\in V$ and $t,s>0$;

(a)\, $\mu(x,t)+\nu(x,t)\leq 1$,

(b)\, $\mu(x,t)>0$,

(c)\, $\mu(x,t)=1 $ if and only if, $x=0$,

(d)\, $\mu(\alpha x,t)=\mu(x,\frac{t}{|\alpha|})$ for each
$\alpha\neq 0$,

(e)\, $\mu(x,t)\ast \mu(y,s)\leq \mu(x+y,t+s)$,

(f)\, $\mu(x,.):(0,\infty)\longrightarrow [0,1]$ is continuous,

(g)\, $lim_{t\longrightarrow\infty}\mu(x,t)=1$ and
$lim_{t\longrightarrow 0}\mu(x,t)=0$,

(h)\, $\nu(x,t)<1$,

(i)\, $\nu(x,t)=0 $ if and only if, $x=0$,

(j)\, $\nu(\alpha x,t)=\nu(x,\frac{t}{|\alpha|})$ for each
$\alpha\neq 0$,

(k)\, $\nu(x,t)\diamond \nu(y,s)\geq \nu(x+y,t+s)$,

(l)\, $\nu(x,.):(0,\infty)\longrightarrow [0,1]$ is continuous,

(m)\, $lim_{t\longrightarrow\infty}\nu(x,t)=0$ and
$lim_{t\longrightarrow 0}\nu(x,t)=1$,
\end{definition} In this case $(\mu,\nu)$
is called a intuitionistic fuzzy norm.
\begin{example} Let $(V,\|.\|)$ be a normed space. Denote $a*b=ab$ and $a\diamond
b=\min(a+b,1)$ for all $a,b\in[0,1]$ and let $\mu_0$ and $\nu_0$
be fuzzy sets on $X^{2}\times (0,\infty)$ defined as follows:
\begin{eqnarray*}\mu_0(x,t)=\frac{t}{t
+\|x\|}&,& \nu_0(x,t)=\frac{\|x\|}{t +\|x\|},\end{eqnarray*} for
all  $t\in \mathbb{R}^+$.
 Then $(V,\mu_0,\nu_0,\ast,\diamond)$ is an intuitionistic fuzzy normed space.
 \end{example}
\begin{definition}A sequence $ \{x_{n}\}$ in a intuitionistic fuzzy normed space $(V,\mu,\nu,\ast,\diamond)$ is called a
Cauchy sequence if for each $\varepsilon >0$ and $ t>0$, there
exists $ n_{0} \in \mathbb{N}$ such that
$$\mu(x_{n}-x_{m},t)>1-\varepsilon,$$ and $$\nu(x_{n}-x_{m},t)<\varepsilon,$$for each $ n,m\geq n_{0}$.
 The sequence $ \{x_{n}\}$ is said to be convergent to $x\in V$
in intuitionistic fuzzy normed space $(V,\mu,\nu,\ast,\diamond)$
and denote by $x_{n}\stackrel{(\mu,\nu)}{\longrightarrow}x$ if  $
\mu(x_{n}-x,t)\longrightarrow 1$ and $
\nu(x_{n}-x,t)\longrightarrow 0$ whenever $n\longrightarrow\infty$
for every $t>0$. A intuitionistic fuzzy normed space is said to be
complete if and only if every Cauchy sequence is
convergent.\end{definition}

\begin{lemma}Let $(V,\mu,\nu,\ast,\diamond)$  be a intuitionistic fuzzy normed space. If we define
$$M(x,y,t)=\mu(x-y,t),$$ and $$N(x,y,t)=\nu(x-y,t),$$ then $(M,N)$ is a intuitionistic fuzzy metric on $V$, which
is induced by the intuitionistic fuzzy norm
$(\mu,\nu)$.\end{lemma}
\begin{lemma}Let $(\mu,\nu)$ be a intuitionistic fuzzy norm, then

(i)\, $\mu(x,t)$ and $\nu(x,t)$ are nondecreasing and
nonincreasing  with respect to $t$, respectively.

(ii)\, $\mu(x-y,t)=\mu(y-x,t)$ and $\nu(x-y,t)=\nu(y-x,t)$ for
every $t>0$.\end{lemma}
\begin{definition}Let $(V,\mu,\nu,\ast,\diamond)$  be a intuitionistic fuzzy normed space. We define open
ball $B(x,r,t)$ with center $x\in V$ and radius $0<r<1$, as
$$B(x,r,t)=\{y\in V: \mu(x-y,t)>1-r, \nu(x-y,t)<r \},t>0.$$ Also a subset
$A\subseteq V$ is called open if for each $x\in A$, there exist
$t>0$ and $0<r<1$ such that $B(x,r,t)\subseteq A$. Let
$\tau_{(\mu,\nu)}$ denote the family of all open subset of $V$.
$\tau_{(\mu,\nu)}$ is called the topology induced by
intuitionistic fuzzy norm.\end{definition} Note that this topology
is the same as the topology induced by intuitionistic fuzzy metric
sense Park (see \cite{pa} Remark 3.3).
\begin{definition}Let $(V,\mu,\nu,\ast,\diamond)$  be a intuitionistic fuzzy normed space. A subset $A$ of $V$ is said
to be IF-bounded if  there exists $t>0$ and $0<r<1$ such that $
\mu(x,t)>1-r$ and $ \nu(x,t)<r$ for each $x\in A$.\end{definition}
\begin{theorem}In a intuitionistic fuzzy normed space every compact set is closed and
IF-bounded.\end{theorem}By Lemma 4.3 the proof is the same as
intuitionistic fuzzy metric spaces (see \cite{pa} Remark 3.10).
\begin{lemma} A subset $A$ of $\mathbb{R}$ is IF-bounded in $(\mathbb{R},\mu,\nu,\ast,\diamond)$ if
and only if is bounded in $\mathbb{R}$.\end{lemma} \textbf{Proof}
Let $A$ is IF-bounded in $(\mathbb{R},\mu,\nu,\ast,\diamond)$,
then there are $t_{0}>0$ and $0<r_{0}<1$ such that for every non
zero $a\in A$ we have $$ 1-r_{0}<\mu(a
,t_{0})=\mu(1,\frac{t_{0}}{|a|}),$$ and
$$ r_{0}>\nu(a
,t_{0})=\nu(1,\frac{t_{0}}{|a|}),$$
 therefore there exists
$k\in \mathbb R^{+}$ such that $|a|\leq k$, that is $A$ is bounded
in $\mathbb R$. The converse is easy.\qed
\begin{lemma} A sequence $\{\beta_{n}\}$ is convergent
 in the intuitionistic fuzzy normed space $(\mathbb{R},\mu,\nu,\ast,\diamond)$ if and only
if it is convergent in $(\mathbb{R},|.|)$ .\end{lemma} \textbf{
Proof} If $|\beta_{n}-\beta|\longrightarrow 0$, then
\begin{eqnarray*}
\lim_{n\longrightarrow\infty}\mu(\beta_{n}-\beta,t)&=&\lim_{n\longrightarrow\infty}\mu(1,\frac{t}{|\beta_{n}-\beta|})\\
& =&\mu(1,\infty)\\&=& 1,\end{eqnarray*} and
\begin{eqnarray*}
\lim_{n\longrightarrow\infty}\nu(\beta_{n}-\beta,t)&=&\lim_{n\longrightarrow\infty}\nu(1,\frac{t}{|\beta_{n}-\beta|})\\
& =&\nu(1,\infty)\\&=& 0,\end{eqnarray*} that is
$\beta_{n}\longrightarrow \beta$ in
$(\mathbb{R},\mu,\nu,\ast,\diamond)$. Conversely suppose that
$\lim_{n\longrightarrow\infty}\mu(\beta_{n}-\beta,t)=1$ and
$\lim_{n\longrightarrow\infty}\nu(\beta_{n}-\beta,t)=0$. If
$\liminf (\beta_{n}-\beta)=u$ and $\limsup (\beta_{n}-\beta)=v$
and $u,v$ are not $+\infty$ or $-\infty$, then we can find
subsequences $\{\beta_{n_{k}}-\beta\}$ and
$\{\beta_{m_{k}}-\beta\}$ converging to $u,v$, respectively. By
assumption, then $\mu(u,t)=\mu(v,t)=1$, for all $t>0$, so $u=v=0$,
i.e. the limit $\{\beta_{n}-\beta\}$ exists and is $0$. If one of
these or both are infinity then since $\mu(x,t)=\mu(1,t/|x|)$ and
$\mu$ is nondecreasing in second variable, then
$$\limsup \mu(1,t/|\beta_{n}-\beta|)\leq
\lim_{n\longrightarrow\infty}\mu(\beta_{n}-\beta,t)\leq \liminf
\mu(1,t/|\beta_{n}-\beta|).$$ Now, if $\liminf
\beta_{n}-\beta=-\infty$ then we have
$$\lim \mu(\beta_{n}-\beta,t)<\liminf \mu(\beta_{n}-\beta,t)=\liminf \mu(1,t/|\beta_{n}-\beta|).$$
This implies $1<0$, by 4.2(g). If $\limsup
\beta_{n}-\beta=+\infty$ then $\liminf \beta-\beta_{n}=-\infty$,
and again $1<0$. Therefore,
$\lim_{n\longrightarrow\infty}(\beta_{n}-\beta)=0$, that is
 $\{\beta_{n}\}$ is
convergent in $(\mathbb R,|.|)$.\qed

By last lemma $(\mathbb{R},\mu,\nu,\ast,\diamond)$ is complete.
\begin{corollary}If the real sequence $\{\beta_{n}\}$ is IF-bounded then
it has at last one limit point.\end{corollary}
\begin{definition} The 5-tuple $(\mathbb R^{n},\Phi,\Psi,\ast,\diamond)$ is called a intuitionistic fuzzy
Euclidean normed space  if $\ast$ is a t-norm, $\diamond$ is a
t-conorm and $(\Phi,\Psi)$ is a intuitionistic fuzzy Euclidean
norm defined by
$$ \Phi(x,t)=\prod_{j=1}^{n}\mu(x_{j},t),$$
and
$$ \Psi(x,t)=\coprod_{j=1}^{n}\nu(x_{j},t),$$
where $x=(x_{1},\cdots,x_{n})$,
$\prod_{j=1}^{n}a_{j}=a_{1}\ast\cdots\ast a_{n}$,
$\coprod_{j=1}^{n}a_{j}=a_{1}\diamond\cdots\diamond a_{n}$, $t>0$,
and $(\mu,\nu)$ is a intuitionistic fuzzy norm  .\end{definition}
\begin{lemma} If $\diamond=\max$ then $(\mathbb R^{n},\Phi,\Psi,\ast,\diamond)$
is an intuitionistic fuzzy  normed space.\end{lemma} We omit the
proof because it is straightforward.
\begin{corollary} The intuitionistic fuzzy Euclidean normed space $(\mathbb R^{n},\Phi,\Psi,\ast,\max)$
is complete.\end{corollary}

\section{FINITE DIMENSIONAL INTUITIONISTIC FUZZY NORMED SPACES}
\begin{theorem}Let $\{x_{1},..., x_{n}\}$ be a linearly independent set
of vectors in vector space $V$ and $(V,\mu,\nu,\ast,\diamond)$  be
a intuitionistic fuzzy normed space. Then there are  numbers
$c,d\neq 0$ and a intuitionistic fuzzy norm space
$(\mathbb{R},\mu_0,\nu_0,\ast,\diamond)$ such that for every
choice of real scalars $\alpha_{1},...,\alpha_{n}$ we have
\begin{equation}\mu(\alpha_{1}x_{1}+...+\alpha_{n}x_{n},t)\leq
\mu_{0}(c[|\alpha_{1}|+...+|\alpha_{n}|],t),\end{equation} and
\begin{equation}\nu(\alpha_{1}x_{1}+...+\alpha_{n}x_{n},t)\geq
\nu_{0}(d[|\alpha_{1}|+...+|\alpha_{n}|],t).\end{equation}
\end{theorem}
\textbf{Proof.} Put $s=|\alpha_{1}|+...+|\alpha_{n}|$. If $s=0$ ,
all $\alpha_{j}$'s must be zero, so (5.1) and (5.2) holds for any
$c,d\neq 0$. Let $s>0$. Then (5.1) and (5.2) are equivalent to the
inequalities which we obtain from (5.1) and (5.2) by dividing by
$s$ and putting $\beta_{j}=\frac{\alpha_{j}}{s}$, that is,
\begin{equation}\mu(\beta_{1}x_{1}+...+\beta_{n}x_{n},t')\leq
\mu_{0}(c,t'),t'=\frac{t}{s},\sum_{j=1}^{n}|\beta_{j}|=1,\end{equation}
and
\begin{equation}\nu(\beta_{1}x_{1}+...+\beta_{n}x_{n},t')\geq
\nu_{0}(d,t'),t'=\frac{t}{s},\sum_{j=1}^{n}|\beta_{j}|=1,\end{equation}
Hence it suffices to prove the existence of a $c,d\neq 0$ and
intuitionistic fuzzy norm $(\mu_{0},\nu_0)$ such that (5.3) and
(5.4) holds. Suppose that this is not true.
 Then there exists a sequence $ \{y_{m}\}$ of vectors,
$$y_{m}=\beta_{1,m}x_{1}+...+\beta_{n,m}x_{n},(\sum_{j=1}^{n}|\beta_{j,m}|=1),$$
such that $\mu(y_{m},t)\longrightarrow 1$ and
$\nu(y_{m},t)\longrightarrow 0$ as $m\longrightarrow\infty$ for
every $t>0$. Since $\sum_{j=1}^{n}|\beta_{j,m}|=1$, we have,
$|\beta_{j,m}|\leq 1$ then by 4.9 the sequence of
$\{\beta_{j,m}\}$ is IF-bounded. In according 4.11
$\{\beta_{1,m}\}$ has a convergent subsequence. Let $\beta_{1}$
denote the limit of that subsequence, and let $\{y_{1,m}\}$ denote
the corresponding subsequence of $\{y_{m}\}$. By the same
argument, $\{y_{1,m}\}$ has a subsequence $\{y_{2,m}\}$ for which
the corresponding of real scalars $\beta_{2,m}$ convergence ; let
$\beta_{2}$ denote the limit.Continuing this process, after $n$
steps we obtain a subsequence $\{y_{n,m}\}_{m}$ of $\{y_{m}\}$
such that$$
y_{n,m}=\sum_{j=1}^{n}\gamma_{j,m}x_{j},(\sum_{j=1}^{n}|\gamma_{j,m}|=1)$$
and $\gamma_{j,m}\longrightarrow\beta_{j}$ as
$m\longrightarrow\infty$.
Since,\begin{eqnarray*}lim_{m}\mu(y_{n,m}-\sum_{j=1}^{n}\beta_{j}x_{j},t)&=&lim_{m}\mu(\sum_{j=1}^{n}
(\gamma_{j,m}-\beta_{j})x_{j},t)\\
& \geq & lim_{m}[\mu((\gamma_{1,m}-\beta_{1})x_{1},t/n)\ast
...\ast \mu((\gamma_{n,m}-\beta_{n})x_{n},t/n)]\\&=&
1,\end{eqnarray*} and
\begin{eqnarray*}lim_{m}\nu(y_{n,m}-\sum_{j=1}^{n}\beta_{j}x_{j},t)&=&lim_{m}\nu(\sum_{j=1}^{n}
(\gamma_{j,m}-\beta_{j})x_{j},t)\\
& \leq & lim_{m}[\nu((\gamma_{1,m}-\beta_{1})x_{1},t/n)\diamond
...\diamond \mu((\gamma_{n,m}-\beta_{n})x_{n},t/n)]\\&=&
0.\end{eqnarray*}

Hence,
$$lim_{m\longrightarrow\infty}y_{n,m}=\sum_{j=1}^{n}\beta_{j}x_{j},\sum_{j=1}^{n}|\beta_{j}|=1,$$
  so that not all $\beta_{j}$ can be zero. Put $y=\sum_{j=1}^{n}\beta_{j}x_{j}$. Since
$\{x_{1},...,x_{n}\}$ is a linearly independent set, we thus have
$y\neq 0$. Since $ \mu(y_{m},t)\longrightarrow 1$ and $
\nu(y_{m},t)\longrightarrow 0$ by assumption then we have $
\mu(y_{n,m},t)\longrightarrow 1$ and $
\nu(y_{n,m},t)\longrightarrow 0$.
Hence,\begin{eqnarray*}\mu(y,t)&=&\mu((y-y_{n,m})+y_{n,m},t)\\
&\geq & \mu(y-y_{n,m},t/2)\ast \mu(y_{n,m},t/2)\\&\longrightarrow
1\end{eqnarray*} and \begin{eqnarray*}\nu(y,t)&=&\nu((y-y_{n,m})+y_{n,m},t)\\
&\leq & \mu(y-y_{n,m},t/2)\diamond
\nu(y_{n,m},t/2)\\&\longrightarrow 0\end{eqnarray*}and so , $y=0$
which is a contradicts.\qed

\begin{definition}Let $(V,\mu,\nu,\ast,\diamond)$ and $(V,\mu',\nu',\ast',\diamond')$ be intuitionistic
fuzzy normed space. Then two intuitionistic  fuzzy norms
$(\mu,\nu)$ and $(\mu',\nu')$ are said to be equivalent whenever
$x_{n}\stackrel{(\mu,\nu)}{\longrightarrow}x$ in
$(V,\mu,\nu,\ast,\diamond)$ If and only if
$x_{n}\stackrel{(\mu',\nu')}{\longrightarrow}x$ in
 $(V,\mu',\nu',\ast',\diamond')$.\end{definition}
\begin{theorem} On a finite dimensional vector space $V$, every two intuitionistic
fuzzy norms $(\mu,\nu)$ and $(\mu',\nu')$  are
equivalent.\end{theorem} \textbf{Proof.} Let $\text{dim} V=n$ and
$\{v_{1},...,v_{n}\}$ be a basis for $V$. Then every $x\in V$ has
a unique representation $ x=\sum_{j=1}^{n}\alpha_{j}v_{j}$. Let
$x_{m}\stackrel{(\mu,\nu)}{\longrightarrow}x$ in
$(V,\mu,\nu,\ast,\diamond)$ but for each $m\in \mathbb{N}$,
$x_{m}$ has a unique representation i.e.
$$x_{m}=\alpha_{1,m}v_{1}+...+\alpha_{n,m}v_{n}.$$
By Theorem 5.1 there are $c,d\neq 0$ and a intuitionistic fuzzy
norm $(\mu_0,\nu_0)$ such that (5.1) and (5.2) hold. So
\begin{eqnarray*}
\mu(x_{m}-x,t)& \leq &
\mu_{0}(c\sum_{j=1}^{n}|\alpha_{j,m}-\alpha_{j}|,t)\\
 &\leq & \mu_{0}(c|\alpha_{j,m}-\alpha_{j}|,t),\end{eqnarray*}and\begin{eqnarray*}
\nu(x_{m}-x,t)& \geq &
\nu_{0}(d\sum_{j=1}^{n}|\alpha_{j,m}-\alpha_{j}|,t)\\
 &\geq & \nu_{0}(d|\alpha_{j,m}-\alpha_{j}|,t).\end{eqnarray*}
 Now if $m\longrightarrow \infty$ then $\mu(x_{m}-x,t)\longrightarrow
  1$ and $\nu(x_{m}-x,t)\longrightarrow
  0$ for every $t>0$ hence
 $|\alpha_{j,m}-\alpha_{j}|\longrightarrow 0$ in $\mathbb{R}$. On the
 other hand,
 \begin{eqnarray*}\mu'(x_{m}-x,t)&\geq &
 \mu'((\alpha_{1,m}-\alpha_{1})v_{1},t/n)\ast' ...\ast'
 \mu'((\alpha_{n,m}-\alpha_{n})v_{n},t/n)\\
 &=& \mu'(v_{1},\frac{t}{n(\alpha_{1,m}-\alpha_{1})})\ast' ...
 \ast'
 \mu'(v_{n},\frac{t}{n(\alpha_{n,m}-\alpha_{n})}).\end{eqnarray*}
 and \begin{eqnarray*}\nu'(x_{m}-x,t)&\leq &
 \nu'((\alpha_{1,m}-\alpha_{1})v_{1},t/n)\diamond' ...\diamond'
 \nu'((\alpha_{n,m}-\alpha_{n})v_{n},t/n)\\
 &=& \nu'(v_{1},\frac{t}{n(\alpha_{1,m}-\alpha_{1})})\diamond' ...
 \diamond'
 \nu'(v_{n},\frac{t}{n(\alpha_{n,m}-\alpha_{n})}).\end{eqnarray*} Since
 $|\alpha_{j,m}-\alpha_{j}|\longrightarrow 0$ so
 $\frac{t}{n(\alpha_{j,m}-\alpha_{j})}\longrightarrow \infty$ so
 $\mu'(v_{j},\frac{t}{n(\alpha_{j,m}-\alpha_{j})})\longrightarrow 1$ and
 $\nu'(v_{j},\frac{t}{n(\alpha_{j,m}-\alpha_{j})})\longrightarrow 0$.
 Then $x_{m}\stackrel{(\mu',\nu')}{\longrightarrow}x$ in $(V,\mu',\nu',\ast',\diamond')$.
 With the same argument $x_{m}\longrightarrow x$ in $(V,\mu',\nu',\ast',\diamond')$
 imply $x_{m}\longrightarrow x$ in $(V,\mu,\nu,\ast,\diamond)$.\qed
\section{BOUNDED LINEAR OPERATORS}
\begin{definition} A linear operator $T:(V,\mu,\nu,\ast,\diamond)\longrightarrow
(V',\mu',\nu',\ast',\diamond')$ is said to be intuitionistic fuzzy
bounded if there exists  constants $h,k\in \mathbb{R}-\{0\}$ such
that for every $x\in V$ and for every $t>0$ ,$$\mu'(Tx,t)\geq
\mu(hx,t),$$ and $$\nu'(Tx,t)\leq \nu(kx,t).$$ \end{definition}
\begin{corollary} Every fuzzy bounded linear operator is
continuous.\end{corollary}
\begin{definition} A linear operator $T:(V,\mu,\nu,\ast,\diamond)\longrightarrow
(V',\mu',\nu',\ast',\diamond')$ is a fuzzy topological isomorphism
if $T$ is one-to-one and onto, and both $T$ and $T^{-1}$ are
continuous. Intuitionistic fuzzy normed spaces
$(V,\mu,\nu,\ast,\diamond)$ and $(V',\mu',\nu',\ast',\diamond')$
for which such a $T$ exists are intuitionistic fuzzy topologically
isomorphic.\end{definition}
\begin{lemma} A linear operator $T:(V,\mu,\nu,\ast,\diamond)\longrightarrow
(V',\mu',\nu',\ast',\diamond')$ is intuitionistic fuzzy
topological isomorphism if $T$ is onto and there exists constants
$a,b,a',b'\neq 0$ such that $\mu(ax,t)\leq \mu'(Tx,t)\leq
\mu(bx,t)$ and $\nu(a'x,t)\leq \nu'(Tx,t)\leq
\nu(b'x,t)$.\end{lemma} \textbf{Proof.} By hypothesis $T$ is
intuitionistic fuzzy bounded and by last corollary is continuous
and since $Tx=0$ implies $1=\mu'(Tx,t)\leq \mu(x,t/|b|)$ and
consequently $x=0$ then $T$ is one-to-one. Thus $T^{-1}$ exists
and, since $\mu'(Tx,t)\leq \mu(bx,t)$  and $\nu'(Tx,t)\leq
\nu(b'x,t)$ are equivalent to $\mu'(y,t)\leq
\mu(bT^{-1}y,t)=\mu(T^{-1}y,t/|b|)$ and $\nu'(y,t)\leq
\nu(b'T^{-1}y,t)=\nu(T^{-1}y,t/|b'|)$  or
$\mu'(\frac{1}{b}y,t)\leq \mu(T^{-1}y,t)$ and
$\nu'(\frac{1}{b'}y,t)\leq \nu(T^{-1}y,t)$ where $y=Tx$, we see
$T^{-1}$ is intuitionistic fuzzy bounded and by last corollary is
continuous. Hence $T$ is a intuitionistic fuzzy topological
isomorphism. \qed
\begin{corollary} Fuzzy topologically isomorphism preserves
completeness.\end{corollary}
\begin{theorem} Every linear operator  $T:(V,\mu,\nu,\ast,\diamond)\longrightarrow
(V',\mu',\nu',\ast, \diamond)$ where $\diamond=\max$ and
$\text{dim}V<\infty$ but other, not necessarily finite
dimensional, is continuous.\end{theorem}\textbf{ Proof.} If we
define
\begin{eqnarray} \mu''(x,t)=\mu(x,t)\ast \mu'(Tx,t),\end{eqnarray} and
\begin{eqnarray}\nu''(x,t)=\nu(x,t)\diamond \nu'(Tx,t),\end{eqnarray}
 then

$(V,\mu'',\nu'',\ast,\diamond)$ is intuitionistic fuzzy normed
spaces because, (a),(b),(c),(d),(f),(g),(h),(l) and (m) immediate
of definition, for triangle inequalities (e) and (k),
\begin{eqnarray*} \mu''(x,t)\ast \mu''(z,s)&=& [\mu(x,t)\ast \mu'(Tx,t)]\ast[\mu(z,s)\ast
\mu'(Tz,s)]\\&=& [\mu(x,t)\ast \mu(z,s)]\ast[\mu'(Tx,t)\ast \mu'(Tz,s)]\\
&\leq & \mu(x+z,t+s)\ast \mu'(T(x+z),t+s)\\
&=& \mu''(x+z,t+s).\end{eqnarray*} The proof for (k) is similar
above. Now, let $x_{n}\stackrel{(\mu,\nu)}{\longrightarrow}x$ then
by Theorem 3.3 $x_{n}\stackrel{(\mu'',\nu'')}{\longrightarrow}x$
but since by (6.1) and (6.2), $\mu'(Tx,t)\geq \mu''(x,t)$ and
$\nu'(Tx,t)\leq \nu''(x,t)$ then
$Tx_{n}\stackrel{(\mu',\nu')}{\longrightarrow}Tx$. Hence $T$ is
continuous.\qed
\begin{corollary} Every linear isomorphism between finite dimensional
fuzzy normed spaces is topological isomorphism.\end{corollary}
\begin{corollary} Every finite dimensional intuitionistic fuzzy normed space
$(V,\mu,\nu,\ast,\diamond)$ where $\diamond=\max$ is
complete.\end{corollary} \textbf{Proof.} Let $\diamond=\max$. By
last corollary $(V,\mu,\nu,\ast,\diamond)$ is topologically
isomorphism to $(\mathbb R^{n},\mu,\nu,\ast,\diamond)$. Since
$(\mathbb R^{n},\mu,\nu,\ast,\diamond)$ is complete and
topological isomorphism preserves completeness then
$(V,\mu,\nu,\ast,\diamond)$ is complete.\qed

\bibliographystyle{amsplain}

\begin{thebibliography}{99}
\bibitem{as1}M. Amini and R. Saadati, Topics in fuzzy metric space, \textit{J. Fuzzy
Math.}, \textbf{4 }(2003) 765-768.
\bibitem{as2}M. Amini and R. Saadati, Some properties of Continuous t-norms and s-norms, \textit{Int.
J. Pure Appl. Math.}, \textbf{16 }(2004)157-164.
\bibitem{ch} Hu, Chengming, $\mathcal{C}$-structure of FTS. V: Fuzzy metric
spaces, \textit{J. Fuzzy Math.},\textbf{ 3}, (1995) 711-721.
\bibitem{el1}M.S. Elnaschie, On the uncertainty of Cantorian geometry and two-slit experiment.
\textsl{Chaos, Soliton and Fractals}, \textbf{9}(1998)517-529.
\bibitem{el2}M.S. Elnaschie,On the verifications of heterotic strings theory and $\epsilon^{(\infty)}$ theory.
\textsl{Chaos, Soliton and Fractals} ,\textbf{11} (2000)397-407.
\bibitem{gv1}A. George and P. Veeramani, On some result in fuzzy metric space, \textit{Fuzzy
Sets and System},  \textbf{64} (1994)395-399.
\bibitem{gv2}A. George and P. Veeramani, On some result of analysis for fuzzy metric
spaces, \textit{Fuzzy Sets and Systems}, \textbf{90}(1997)365-368.
\bibitem{gr1}V. Gregori and S. Romaguera, Some properties of fuzzy
metric spaces, \textit{Fuzzy Sets and Systems}, \textbf{115}
(2000)485-489.
\bibitem{gr2}V. Gregori and S. Romaguera, On completion of fuzzy
metric spaces, \textit{Fuzzy Sets and Systems}, \textbf{130}
(2002)399-404.
\bibitem{gr3}V. Gregori, S. Romaguera, Characterizing completable
fuzzy metric spaces, \textit{Fuzzy Sets and Systems}, \textbf{144}
(2004) 411-420.
\bibitem{j}K.D. Joshi, Introduction to General Topology, Wiely
Estern, Bombay, 1991.
\bibitem{k}K. Kuratowski, Topology,Academic press, New York, 1966.
\bibitem{lo}R. Lowen, Fuzzy set theory, Kluwer Academic Publishers,  Dordrecht, 1996.
\bibitem{m}R.E. Megginson, An Introduction to Banach space theory,
Springer-Verlag, New York, 1998.
\bibitem{pa} J. H. Park, Intuitionistic fuzzy metric spaces, \textit{Chaos, Solitons and Fractals} \textbf{22}
(2004) 1039-1046.
\bibitem{lr} J. Rodr\'{i}guez-L\'{o}pez, S. Ramaguera, The Hausdorff fuzzy metric on compact
sets, \textit{Fuzzy Sets and Systems}, in press.
\bibitem{sa} R. Saadati and S.M. Vaezpour, Some results on fuzzy Banach
spaces, to appear in  \emph{J. Appl. Math. Comput}.
\bibitem{sw}B. Schweizer and A. Sklar, Satistical metric spaces,
\textit{Pacific J. Math.},\textbf{10} (1960)314-334.
\bibitem{za}L.A. Zadeh, Fuzzy sets, \textit{Inform. and control}, \textbf{8} (1965)338-353.
\end{thebibliography}

\end{document}